\documentclass[11pt,psamsfonts]{amsart}
\usepackage{amsmath}
\usepackage{amsthm}
\usepackage{amssymb}
\usepackage{amscd}
\usepackage{amsfonts}
\usepackage{amsbsy}
\usepackage{epsfig,afterpage}
\usepackage{psfrag}
\usepackage{graphicx}
\usepackage{indentfirst, latexsym, bm,amssymb}
\usepackage{bbding}

\def\v{\varphi}

\newtheorem {theorem} {Theorem}

\newcommand{\R}{\mathbb{R}}

\newcommand{\X}{{\mathcal X}}

\begin{document}

\title[Global structure of quaternion differential equations]
{Global structure of quaternion  polynomial differential
equations}

\author[Xiang Zhang]
{ Xiang Zhang }\footnote{The author is partially supported by NNSF
of China grant 10831003 and Shanghai Pujiang Program grant
09PJD013.}

\address{ Department of Mathematics, Shanghai Jiaotong University, Shanghai
200240, P. R. China} \email{xzhang@sjtu.edu.cn}



\date{}

\begin{abstract}
In this paper we mainly study the global structure of the
quaternion Bernoulli equations $\dot q=aq+bq^n$ for $q\in \mathbb
H$ the quaternion field and also some other form of cubic
quaternion differential equations. By using the Liouvillian
theorem of integrability and the topological characterization of
$2$--dimensional torus: orientable compact connected surface of
genus one, we prove that the quaternion Bernoulli equations may
have invariant tori, which possesses a full Lebesgue measure subset of
$\mathbb H$. Moreover, if $n=2$ all the invariant tori are full of
periodic orbits; if $n=3$ there are infinitely many invariant tori
fulfilling periodic orbits and also infinitely many invariant ones
fulfilling dense orbits.

\end{abstract}

\maketitle

\smallskip

\section{Introduction and main results}\label{s1}

\smallskip

\smallskip

The dynamics of ordinary differential equations in $\mathbb R$ or
$\mathbb C$ has been intensively studied from many different
points of view. While because of the noncommutativity of the
quaternion algebra, the study on quaternion differential equations
becomes very difficult and much involved, and the results in this
field are very few. Recent years because of their application in
quantum and fluid mechanics, see e.g. \cite{Ad1, Ad2, FJSS, Gi,
GHKR, RR01, RR97},  the study on the dynamics of quaternion
differential equations has been attracting more interesting.

\smallskip

In 2006 Campos and Mawhin \cite{CM} initiated the study on the
existence of periodic solutions of one--dimensional first order
periodic quaternion differential equations. Wilczy\'nski
\cite{Wi09} continued this study and payed more attention on the
existence of two periodic solutions of quaternion Riccati
equations. Our work in \cite{GLZ} presented a study on the global
structure of the quaternion autonomous homogeneous differential
equations
\begin{equation}\label{e1}
\dot q=aq^n, \qquad q\in\mathbb H,
\end{equation}
where $a\in\mathbb H$ is a parameter. Recall that $\mathbb{H}$ is
the quaternion field.

\smallskip

In this paper we will study the global dynamics of the quaternion
Bernoulli equations
\begin{equation}\label{e2}
\dot q=bq+aq^n,
\end{equation}
with $a,b,q\in\mathbb H$, $2\le n\in\mathbb N$ and also of the
third order equation
\begin{equation}\label{e3}
\dot q=a(q-c_0)(q+c_0)q,
\end{equation}
with $a\in\mathbb H$ and $c_0\in \mathbb R$.

\smallskip

The quaternion Bernoulli equation \eqref{e2} consists of linear
terms and homogenous nonlinearities of degree $n$. We note that
real planar polynomial vector fields generalizing the linear
systems with homogeneous nonlinearities have been extensively
studied from different points of view, for instance limit cycles,
centers, phase portraits and integrability, see e.g.
\cite{GL04,GLMM,LLLZ,LV09}. Some famous three dimensional real
differential systems exhibiting chaotic phenomena, for instance
Lorenz system, Rabinovich systems and Rikitake systems and so on,
also have this form, which consist of linear terms and homogeneous
nonlinearity of degree 2. As our knowledge the dynamics of the
quaternion equations of form \eqref{e2} with $a, b\ne 0$ has never
been studied. We note that for either $a=0$ or $b=0$, equation
\eqref{e2} is in fact the equation \eqref{e1}, and it has been
studied in \cite{GLZ}.

\smallskip

Equation \eqref{e2} with $a,b\ne 0$ can be written in
\begin{equation}\label{e5}
\dot q=a(cq-q^n),
\end{equation}
with $a,c\in\mathbb H$ not zero.

\smallskip

Our first result is the following.

\smallskip

\begin{theorem}\label{t1} For the quaternion differential equation \eqref{e5} with
$c\in \mathbb R$ not zero, the following statements hold.
\begin{itemize}

\smallskip

\item[$(a)$] Assume that $a+\overline a\ne 0$ and $a-\overline
a=0$.

\smallskip

\begin{itemize}

\item[$(a_1)$] The phase space $\mathbb R^4$, i.e., $\mathbb H$,
is foliated by invariant planes of \eqref{e5}, which all pass
through the origin.

\smallskip

\item[$(a_2)$] On each invariant plane, there are $n$
singularities: one is the origin and the others are located on the
circle centered at the origin with the radius $\sqrt[n-1]{|c|}$,
denoted by $S_c$. All non--trivial orbits are heteroclinic, and
connect the origin and one of the singularities on $S_c$ except
the following $2(n-1)$ ones: there are exactly $n-1$ heteroclinic
orbits connecting the origin and the infinity, and also $n-1$ ones
connecting each one of the singularities on $S_c$ and the
infinity.

\smallskip

\end{itemize}

\item[$(b)$] Assume that $a^2-{\overline a}^2\ne 0$.

\begin{itemize}

\smallskip

\item[$(b_1)$] Each orbit of system \eqref{e5} starting on the
branch of $P:=\{q^{n-1}+\overline q^{n-1}-c=0\}$ is heteroclinic
connecting the origin and one of the singularities given by
$q^{n-1}=c$, which are located in two consecutive region limited
by the branches of $P$.

\smallskip

\item[$(b_2)$] There exists at least one orbit in each connected
region limited by the branches of $P$, which connects the infinity
and one of the singularities of \eqref{e5}.

\end{itemize}

\smallskip

\item[$(c)$] Assume that $a+\overline a=0$ and $a-\overline a\ne
0$.

\smallskip

\begin{itemize}

\item[$(c_1)$] The hypersurfaces $P$ are invariant, on which all
orbits are nontrivial and located in two dimensional invariant
algebraic varieties.

\smallskip

\item[$(c_2)$] The invariant set $\R^4\setminus\{P\}$ is foliated
by one invariant plane, two $2$--dimensional invariant algebraic  varieties
and $2$--dimensional invariant tori. The invariant plane is
foliated by $n$ isochronous centers with $n$ separatrices going to
infinity. One of the algebraic varieties is full of singularities and the
other fulfils periodic orbits with a center and finitely many
heteroclinic orbits.
\end{itemize}
\end{itemize}
\end{theorem}

\smallskip

A {\it nontrivial orbit} is an orbit which is not a singularity. An
{\it algebraic variety} is a subset of $\mathbb R^4$ formed by the common
zeros of finitely many polynomials.

\smallskip

In statement $(c_2)$ of Theorem \ref{t1} we do not study the
dynamics of equation \eqref{e5} on the invariant tori. In fact,
the next theorem shows that the dynamics on the invariant tori
depend on the degree $n$ of the equations.

\smallskip

Now we study the dynamics of equation \eqref{e5} on the invariant
tori appearing in statement $(c_2)$ of the last theorem for
$n=2,3$. For larger $n$, we have no methods to tackle it. The
difficulty is the parametrization of the invariant tori as we will
see in the proof of the following results.

\begin{theorem}\label{t1.1} For the $2$--dimensional invariant tori stated
in $(c_2)$ of Theorem $\ref{t1}$ the following statements hold.

\begin{itemize}

\item[$(a)$]  $n=2$. Each torus is full of periodic orbits.

\smallskip

\item[$(b)$]  $n=3$. Among the tori there are infinite many ones
fulfilled periodic orbits and also infinite many ones fulfilled
dense orbits.
\end{itemize}
\end{theorem}

\smallskip

The above results are on equation \eqref{e5} with $c\in\mathbb R$.
We  now study the equation with $c\in\mathbb H\setminus\mathbb R$.
For general $a\in\mathbb H$ and $2<n\in\mathbb N$, we have no
method to deal with it. The next result is on equation \eqref{e5}
with $0\ne a\in\mathbb R$ and $n=2$.

\smallskip

\begin{theorem}\label{th2.5} For equations \eqref{e5} with
$a\in\mathbb R$ nonzero, $n=2$ and $c-\overline c\ne 0$, set
$L=c_0 q_0+c_1q_1+c_2q_2+c_3q_3-(c_0^2+c_1^2+c_2^2+c_3^2)/2$, the
following statements hold.
\begin{itemize}

\smallskip

\item[$(a)$] If $c+\overline c\ne 0$, all the orbits of system
\eqref{e5} starting on the hyperplane $L=0$ are heteroclinic and
spirally approach the singularities $O=(0,0,0,0)$ and
$S=(c_0,c_1,c_2,c_3)$. There are other two heteroclinic orbits
connecting the infinity and either $S$ or $O$.

\smallskip

\item[$(b)$] If $c+\overline c=0$, the hyperplane $L=0$ is
invariant. The invariant set $\R^4\setminus\{L=0\}$ is foliated by
one invariant plane foliated by two period annuli, one invariant
sphere fulfilling periodic orbits, and $2$--dimensional invariant
tori.
\end{itemize}
\end{theorem}

\smallskip

We remark that the case $c-\overline c=0$ was studied in Theorems
\ref{t1} and \ref{t1.1}.

\smallskip

Finally we study the cubic quaternion differential equation
\eqref{e3}. Without loss of generality we assume $c_0> 0$.

\smallskip

\begin{theorem}\label{th2.6}
Consider the cubic equation \eqref{e3} with $c_0>0$ and
$a\in\mathbb H$ nonzero. Set $L=q_0^2-q_1^2-q_2^2-q_3^2-c_0^2/2$
and denote by $L^+$ and $L^-$ the two sheets of the generalized
hyperboloid of $L=0$ corresponding to $q_0\ge c_0/\sqrt 2$ and
$q_0\le -c_0/\sqrt 2$, respectively. The following statements
hold.
\begin{itemize}

\smallskip

\item[$(a)$] If $a+\overline a\ne 0$, any orbit starting on $L^+$
$($resp. $L^-)$ is heteroclinic connecting the singularities
$O=(0,0,0,0)$ and $S_+=(c_0,0,0,0)$ $($resp. $O$ and
$S_-=(-c_0,0,0,0))$.

\smallskip

\item[$(b)$] If $a+\overline a= 0$, the hyperboloid $L=0$ is
invariant under the flow of \eqref{e3}.  The invariant subset
$\R^4\setminus \{L=0\}$ is foliated by periodic orbits and
$2$--dimensional invariant tori. Of the invariant tori, there are
infinitely many ones fulfilling periodic orbits and also
infinitely many ones fulfilling dense orbits.
\end{itemize}
\end{theorem}

\smallskip

From Theorems \ref{t1.1} and \ref{th2.6} we conjecture that {\it
for quaternion polynomial differential equations of degree larger
than $2$, if the equations have invariant tori, then of which there
are infinitely many ones fulfilling periodic orbits and also
infinitely many ones fulfilling dense orbits}.

\smallskip

We remark that in the proof of the existence of invariant tori, we
will use both the Liouvillian theorem of integrability and also
the topological characterization of torus. In the case that the
mentioned equations have two functionally independent first
integrals but they are not Liouvillian integrable, we prove the
existence of invariant tori by showing that the connected parts of
the intersection of the level sets of the two first integrals are
orientable compact smooth surfaces of genus one.

\smallskip

In this paper, as a by product of our resusts we find some new class of integrable systems. The problem on searching
integrable differential equations, including the integrable Hamiltonian systems, has a long history. It can be traced back
to Poincar\'e and Darboux, and even earlier. In recent years Calogero has done a series of researches
in this direction, see for instance \cite{BC04, CD06, CD04, IC02} and the reference therein.

\smallskip

The paper is organized as follows. In the next section we recall
some basic facts on quaternion which will be used later on. In
Section \ref{s4} we will prove our main results. The last section
is the appendix presenting the results for linear quaternion
equations.

\bigskip

\section{Basic preliminaries}\label{s0}

In this section for readers' convenience we recall some basic
facts on quaternion  algebra (see e.g., \cite{Fr, Ha, Han}), which
will be used later on. Quaternions are non--commutative extension
of complex numbers, which are defined as the field
\[
\mathbb H=\{q=q_0+q_1i+q_2j+q_3k;\,q_0,q_1,q_2,q_3\in\R\},
\]
with $i,j,k$ satisfying
\[
i^2=j^2=k^2=-1, \quad ij=-ji=k.
\]
For $a,b\in\mathbb H$, their addition and multiplication are
defined respectively as
\begin{eqnarray*}
a+b&=&(a_0+b_0)+(a_1+b_1)i+(a_2+b_2)j+(a_3+b_3)k,\\
ab&=&(a_0b_0-a_1b_1-a_2b_2-a_3b_3)
+(a_1b_0+a_0b_1-a_3b_2+a_2b_3)i\\
&&+(a_2b_0+a_3b_1+a_0b_2-a_1b_3)j+(a_3b_0-a_2b_1+a_1b_2+a_0b_3)k.
\end{eqnarray*}
Obviously $a,b\in\mathbb H$ commute if and only if the
vectors $(a_1,a_2,a_3)$ and $(b_1,b_2,b_3)$ are parallel in
$\mathbb R^3$.

\smallskip

For $a\in\mathbb H$, its conjugate is $\overline
a=a_0-a_1i-a_2j-a_3k$. Then we have $\overline {ab}=\overline
b\,\overline a$, $ab+\overline b\,\overline a=ba+\overline a\,
\overline b$ and $a\overline a=a_0^2+a_1^2+a_2^2+a_3^2$. The last
equality implies that $(a_1i+a_2j+a_3k)^2=-(a_1^2+a_2^2+a_3^2)$.

\smallskip

For any $a\in\mathbb H$ nonzero, $\overline a/(a\overline a)$ is
its unique inverse, denoted by $a^{-1}$. Moreover, it is easy to
check that the elements in $\mathbb H$ satisfy the law of
association and distribution under the action of the addition and
multiplication.

\smallskip

Mostly we will use the quaternion structures to prove our results.
But sometimes it is not enough in the proof, we need to write the
quaternion differential equations in components. Considering
one--dimensional quaternion ordinary differential equations
\begin{equation}\label{e1.1}
\dot q=\frac{dq}{dt}=f(q,\overline q),\qquad q\in\mathbb H,
\end{equation}
where $f(q,\overline q)$ is an  $\mathbb H$--valued function in
the variables $q$ and $\overline q$. Set
\[
f(q,\overline q)=f_0(q^*)+f_1(q^*)i+f_2(q^*)j +f_3(q^*)k,
\]
where $q=q_0+q_1i+q_2j+q_3k$ and $q^*=(q_0,q_1,q_2,q_3)\in\R^4$.
Then equation \eqref{e1.1} can be written in an equivalent way as
\[
\dot q_s=f_s(q_0,q_1,q_2,q_3) \qquad\mbox{ for } s=0,1,2,3.
\]

\smallskip

Last paragraph shows that a one--dimensional quaternion ordinary
differential equation is in fact equivalent to a system of
four--dimensional real ordinary differential equations. It is well
known that  the dynamics of higher dimensional real differential
systems is usually very difficult to study. Sometimes the
existence of suitable invariants is very useful in the study.
First integral and invariant algebraic hypersurface are two
important invariants. A real valued differentiable function
$H(q,\overline q)$ is a {\it first integral} of \eqref{e1.1} if
the derivative of $H$ with respect to the time $t$ along the
solutions of \eqref{e1.1} is identically zero. An {\it invariant
algebraic hypersurface} of \eqref{e1.1} is defined by the
vanishing set of a real polynomial $F(q,\overline q)$ satisfying
\[
\left.\frac{d F(q, \overline q)}{dt}\right|_{\eqref{e1.1}}=
K(q,\overline q) F(q,\overline q),
\]
with the {\it cofactor} $K(q,\overline q)$ a real polynomial.

\smallskip

In this paper the most difficult part is the search of invariant
algebraic hypersurfaces and of first integrals. Having them we can
obtain the dynamics of the equations with the help of qualitative
methods. This idea can be found in the study of the Lorenz system
\cite{LZ1}, of the Rabinovich system \cite{CCZ} and of the
Einstein-Yang-Mills Equations \cite{LY} and so on.

\bigskip

\section{Proof of the main results}\label{s4}

\smallskip

\subsection{ Proof of Theorem \ref{t1}}

{\it Statement $(a)$}.  Under the assumption of the theorem we
assume without loss of generality that $a=1$, and set
$c=c_0\in\mathbb R$. Then system \eqref{e5} can be written in
\begin{equation}\label{e6}
\dot q=c_0q-\frac{q^n+\overline q^n}{2}-\frac{q^n-\overline
q^n}{q-\overline q}(q_1i+q_2j+q_3k),
\end{equation}
for $q-\overline q\ne 0$, where we have used the fact that
$q-\overline q=2(q_1i+q_2j+q_3k)$ and
\[
q^n=\frac{q^n+\overline q^n}{2}+\frac{q^n-\overline
q^n}{q-\overline q}(q_1i+q_2j+q_3k).
\]
Obviously, $q^n+\overline q^n$ and $(q^n-\overline
q^n)/(q-\overline q)$ are real. Furthermore using the Darboux
theory of integrability we can check easily that
\[
H_2=\frac{q_2}{q_1}, \qquad H_3=\frac{q_3}{q_1},
\]
are two first integrals of equation \eqref{e6}, which follows from
the facts that $q_1=0,q_2=0$ and $q_3=0$ are three invariant
algebraic hyperplanes with the same cofactor $c_0-(q^n-\overline
q^n)/(q-\overline q)$. For more information on the Darboux theory
of integrability, see for instance \cite{Ll,LZ}. 

We remark that
the Darboux theory of integrability was developed for polynomial vector fields in $\mathbb C^n$ and $\mathbb R^n$. Here we can use this
theory in the non--commutative field $\mathbb H$,
because the mentioned invariant algebraic hyperplanes and their cofactors are all real. Generally, if a polynomial differential equation in $\mathbb H$ has its invariant algebraic hypersurfaces all real, we can apply the Darboux theory of integrability by using these hypersurfaces.

\smallskip

The existence of the two functionally independent first integrals
shows that the $\mathbb R^4$ space is foliated by invariant planes
given by $\{H_2=h_2\}\cap\{H_3=h_3\}$ with $h_2,h_3\in \mathbb
R\cup\{\infty\}$. This proves statement $(a_1)$.

\smallskip

We now prove statement $(a_2)$, that is, study the dynamics of
equation \eqref{e6} on each invariant plane.

\smallskip

For any $h_2,h_3\in\mathbb R$, restricted to each invariant plane
$P_{23}:=\{H_2=h_2\}\cap\{H_3=h_3\}$ equation \eqref{e6} becomes
\begin{equation}\label{e7}
\begin{array}{l}
\dot
q_0=c_0q_0-\sum\limits_{s=0}\limits^{[n/2]}\left(\begin{array}{c}
n\\ 2s\end{array}\right) \left(-\Delta^2\right)^sq_0^{n-2s},\\
\dot
q_1=c_0q_1-\sum\limits_{s=1}\limits^{[(n+1)/2]}\left(\begin{array}{c}
n\\ 2s-1\end{array}\right)
\left(-\Delta^2\right)^{s-1}q_0^{n-2s+1},
\end{array}
\end{equation}
where $[\cdot]$ denotes the integer part function,
$\Delta^2=q_1^2+q_2^2+q_3^2=q_1^2(1+h_2^2+h_3^2)$ and we have used
the binormal expansion
\begin{eqnarray*}
q^n&=&\sum\limits_{s=0}\limits^{[n/2]}\left(\begin{array}{c} n\\
2s\end{array}\right)
\left(-\Delta^2\right)^sq_0^{n-2s}\\
&&+\sum\limits_{s=1}\limits^{[(n+1)/2]}\left(\begin{array}{c} n\\
2s-1\end{array}\right)
\left(-\Delta^2\right)^{s-1}q_0^{n-2s+1}(q_1i+q_2j+q_3k),
\end{eqnarray*}
and the fact that $(q_1i+q_2j+q_3k)^2=-\Delta^2$.

\smallskip

For studying the dynamics of equation \eqref{e7} we transfer it to
the complex field. Set $z=q_0+q_1\sqrt{1+h_2^2+h_3^2}\, i$. Then
equation \eqref{e7} can be written in a one dimensional complex
differential equation
\begin{equation}\label{e8}
\dot z=c_0 z-z^n.
\end{equation}
Clearly, this last equation has $n$ singularities in $\mathbb C$:
$z_0=0$ and
$z_k=\sqrt[n-1]{|c_0|}\exp\left(i(\frac{\delta\pi}{n-1}+\frac{2(k-1)\pi}{n-1})\right)$
for $k=1,\ldots,n-1$, where $\delta=0$ if $c_0>0$ or $\delta=1$ if
$c_0<0$. These singularities are all nodes (see e.g. \cite{AGP}),
and $z_0=0$ has different stability than the other $n-1$ ones. By
introducing the polar coordinates $z=re^{i\theta}$ we can prove
that equation \eqref{e8} has exactly
 $n-1$ heteroclinic orbits connecting the origin and the
 infinity, and the unique heteroclinic orbit connecting each $z_k$
 for $k=1,\ldots,n-1$,
 and the infinity. All the other orbits are heteroclinic and
 connect the origin and one of the $z_k's$.
This proves statement $(a_2)$, and consequently statement $(a)$.

\smallskip

For proving statements $(b)$ and $(c)$, we note that for any
$a\in\mathbb H$ there exists a $c\in\mathbb H$ such that
$cac^{-1}=a_0+\sqrt{a_1^2+a_2^2+a_3^2}i$. Moreover equation
\eqref{e5} with $c_0\in \mathbb R$ can be transformed to $\dot
p=cac^{-1}(c_0p-p^n)$ by the change of variables $p=cqc^{-1}$. So
in what follows we assume without loss of generality that
\[
a=a_0+a_1i.
\]

\smallskip

Set
\begin{eqnarray*}
H=\frac{(q\overline q)^{n-1}}{q^{n-1}+\overline
q^{n-1}-c_0},\qquad S= q^{n-1}+\overline q^{n-1}-c_0.
\end{eqnarray*}
We claim that the derivatives of $H$ and $S$ along equation
\eqref{e5} are
\begin{eqnarray}
\left.\frac{dH}{dt}\right|_{\eqref{e5}} &=& (n-1)(a+\overline a)(c_0-H)H,\label{d1}\\
\left.\frac{dS}{dt}\right|_{\eqref{e5}} &=&
(n-1)\left(aq^{n-1}(c_0-q^{n-1})+(c_0-\overline q^{n-1})\overline
q^{n-1}\overline a\right).\label{d2}
\end{eqnarray}
Indeed,
\begin{eqnarray*}
&&\left.\frac{d(q^{n-1}+\overline
q^{n-1})}{dt}\right|_{\eqref{e5}}=\sum\limits_{l=0}\limits^{n-2}\left(q^l\dot
q q^{n-2-l}+\overline q^{n-2-l}\dot{\overline q}\,\overline
q^l\right)\\
&& =\sum\limits_{l=0}\limits^{n-2}\left(q^lac_0q^{n-1-l}+\overline
q^{n-1-l}c_0\overline a\,\overline q^l-q^laq^{2n-2-l}-\overline
q^{2n-2-l}\,\overline a\,\overline q^l \right)\\
&& =(n-1)\left(c_0(aq^{n-1}+\overline q^{n-1}\overline
a)-(aq^{2n-2}+\overline q^{2n-2}\,\overline a)\right).
\end{eqnarray*}
In the last equality we have used the fact that $q^la
q^k+\overline q^k\overline a\,\overline q^l=aq^{k+l}+\overline
q^{k+l}\,\overline a$. This proves equality \eqref{d2}. Using
equality \eqref{d2} and the fact that $q\overline q$ and
$q^{n-1}+\overline q^{n-1}-c_0$ are real, we can prove easily the
equality \eqref{d1}. This proves the claim.

\smallskip

Restricted to the hypersurface $P:=\{S=0\}$ equation \eqref{d2}
becomes
\begin{equation}\label{d3}
\left.\frac{dS}{dt}\right|_{\eqref{e5},P} = (n-1)(a+\overline
a)(q\overline q)^{n-1}.
\end{equation}
For convenience to the following proof, we write equation
\eqref{e5} in a system.
\begin{equation}\label{d4}
\displaystyle\begin{array}{l} \displaystyle\dot
q_0=a_0\left(c_0q_0-\frac{q^n+\overline
q^n}{2}\right)-a_1\left(c_0-\frac{q^n-\overline q^n}{q-\overline
q}\right)q_1,\\
\displaystyle\dot q_1=a_0\left(c_0-\frac{q^n-\overline
q^n}{q-\overline
q}\right)q_1+a_1\left(c_0q_0-\frac{q^n+\overline q^n}{2}\right),\\
\displaystyle\dot q_2=(a_0q_2-a_1q_3)\left(c_0-\frac{q^n-\overline
q^n}{q-\overline
q}\right),\\
\displaystyle\dot q_3=(a_1q_2+a_0q_3)\left(c_0-\frac{q^n-\overline
q^n}{q-\overline q}\right)
\end{array}
\end{equation}

\smallskip

\noindent {\it Statement $(b)$}. By the assumption we can assume
that $a_0=1$. From \eqref{d3} we get that if an orbit of
\eqref{e5} passes through $P$, it should intersect $P$
transversally. So each region limited by the branches of $P$ is
either positively or negatively invariant.

\smallskip

Since $S+c_0$ is a homogeneous polynomial in
$q^*=(q_1,q_2,q_3,q_4)\in\mathbb R^4$, it follows that each branch
of $P$ is either a hyperplane or a generalized hyperboloid. From
the expression of $H$ it follows that each branch of the level
hypersurfaces $H=h$ for $h\in \mathbb R$ (if exist) is compact.

\smallskip

Obviously the level set $H=0$ is the origin, and the level set
$H=c_0$ consists of the roots of $q^{n-1}=c_0$, because $H=c_0$ is
equivalent to $(q^{n-1}-c_0)(\overline q^{n-1}-c_0)=0$. In fact,
these level sets are exactly formed by the singularities. In
addition the compact hypersurfaces $H=h$ approach $P$ when
$h\rightarrow \pm \infty$.

\smallskip

The facts from the last paragraph and equation \eqref{d1} imply
that each orbit starting on $P$ will finally approach two
singularities, which are located in two consecutive regions
limited by $P$. Furthermore, since the function $H$ has different
signs in the two consecutive regions limited by $P$, it follows
from the continuation of $H$ in each region limited by $P$ that
any heteroclinic orbit should go to the level set $H=0$, i.e. the
origin. This proves statement $(b_1)$.

\smallskip

As a by product of the last results we get that for $0<h<c_0$ the
level set $H=h$ is empty.

\smallskip

Statement $(b_2)$ follows from the proof of statement $(b_1)$,
especially the fact that the orbits starting on two consecutive
branches of $P$ are either all get into or all go out the region
limited by the two branches.

\smallskip

\noindent {\it Statement $(c)$}. The assumption means that $a_0=0$
and $a_1\ne 0$. Without loss of generality we take $a_1=1$.

\smallskip

Set
\[
F=q_2^2+q_3^2.
\]
Then $F$ is a first integral of \eqref{e5}, which follows easily
from \eqref{d4} with $a_0=0$. Moreover, we get from \eqref{d1}
that $H$ is also a first integral of \eqref{e5}, which is
functionally independent with $F$. From \eqref{d3} it follows that
each branch of the hypersurface $P$ is invariant.

\smallskip

We first study the dynamics of \eqref{e5} on $P$. For $n=2$  the
level set $P$ is a hyperplane, on which all orbits are parallel
straight lines. For $n>2$ some easy calculations show that
$(q^n-\overline q^n)/(q-\overline q)\ne q^{n-1}+\overline
q^{n-1}$. This verifies that system \eqref{d4} on $P$ has no
singularities. Moreover each orbit on $P$ is located on a cylinder
$F=f>0$ and rotates strictly along the cylinder. This proves
$(c_1)$.

\smallskip

By some direct calculations and using the equality $\overline i
q^n+\overline q^ni=L_{n-1}(\overline iq+\overline q i)$, we get
that
\begin{eqnarray*}
\nabla H&=&\frac{2(n-1)(q\overline q)^{n-2}}{(q^{n-1}+\overline
q^{n-1}-c_0)^2}\times\\
&&\left(\frac{q^n+\overline
q^n}{2}-c_0q_0,(L_{n-1}-c_0)q_1,(L_{n-1}-c_0)q_2,(L_{n-1}-c_0)q_3
\right).
\end{eqnarray*}
So, in the invariant space $\mathbb R^4\setminus\{P\}$, the
critical points of $(H,F)$ form the invariant plane
$S_1:=\{q_2=0\}\cap \{q_3=0\}$, the invariant varieties
$S_2:=\{(q^n+\overline q^n)/2-c_0q_0=0\}\cap\{L_{n-1}=c_0\}$ and
$S_3:=\{(q^n+\overline q^n)/2-c_0q_0=0\}\cap\{q_1=0\}$, where
$L_{n-1}=(q^n-\overline q^n)(q-\overline q)$.

\smallskip

We get from \eqref{d4} that the invariant variety $S_2$ is full of
singularities and that the invariant variety $S_3$ is full of
periodic orbits with a center and the heteroclinic orbits
connecting the singularities on $L_{n-1}=c_0$.

\smallskip

On the invariant plane $S_1$, system \eqref{d4} is simply
\begin{eqnarray*}
\dot
q_0&=&-c_0q_1+\sum\limits_{s=1}\limits^{[(n+1)/2]}\left(\begin{array}{c}n\\2s-1\end{array}\right)
\left(-q_1^2\right)^{s-1}q_0^{n-2s+1},
\\
\dot q_1&=&\,\,\,\,c_0q_0-\sum\limits_{s=0}\limits^{[n/2]}\left(\begin{array}{c}n\\
2s\end{array}\right) \left(-q_1^2\right)^{s}q_0^{n-2s}.
\end{eqnarray*}
Taking $z=q_0+iq_1$, the last equation can be written in
\begin{equation}\label{d44}
\dot z=i\left(c_0z-z^n\right).
\end{equation}
Clearly equation \eqref{d44} has $n$ singularities: one is at the
origin and the others are located on the circle
$|z|=\sqrt[n-1]{|c_0|}$. Applying Theorem 2.1 of \cite{AGP} to
these singularities we get that the origin is an isochronous
center with the period $2\pi/|c_0|$, and the other singularities
are also isochronous centers with the common period
$2\pi/((n-1)|c_0|)$. Furthermore, the periodic orbits surrounding
the origin have different orientation than the ones around the
other singularities. Hence we have obtained the dynamics of
equation \eqref{d44}, and consequently that of equation \eqref{d4} on
the critical sets.

\smallskip

For all regular values $(h,f)$ of $(H,F)$, the intersection
$M_{h,f}=\{H=h\}\cap\{F=f\}$ is a two dimensional compact
invariant manifold, because $M_{h,f}$ does not contain
singularities and the intersection is transversal. We claim that
the connected submanifolds of $M_{h,f}$ are all invariant tori.
Indeed, system \eqref{d4} can be written in a Hamiltonian system
with the Hamiltonian $H$ under the Poisson bracket
$\{\cdot,\cdot\}$ defined by
\[
\{P,Q\}=\nabla P M(q) \nabla Q,
\]
where $P,Q$ are two arbitrary smooth functions in $R^4$ and
\[
M(q)=\frac{(q^{n-1}+\overline q^{n-1}-c_0)^2}{2(n-1)(q\overline
q)^{n-2}}\left(\begin{array}{cccc}0 & -1 & 0 & 0\\ 1 & 0 & 0 &
0\\0 & 0 & 0 & -1\\
0 & 0 & 1 & 0\end{array}\right).
\]
Furthermore the first integrals $H$ and $F$ are in involution
under the Poisson bracket. Then the claim follows from the classic
Liouvillian theorem on integrability. For  more information on
Poisson structures and Liouvillian integrability, see for instance
\cite{AM,Ar}. This proves statement $(c)$.

\smallskip

We complete the proof of the theorem. $\Box$

\smallskip

\subsection{Proof of Theorem \ref{t1.1}}

As in the proof of statement $(c)$ of Theorem \ref{t1}  we take
$a=i$, $c_0>0$ and use the notations given there.

\smallskip

\noindent {\it Statement $(a)$}. Equation \eqref{d4} with $n=2$
becomes
\begin{equation}\label{en3.2}
\begin{array}{ll}
\dot q_0=\,\,\,\,\,(2q_0-c_0)q_1,\qquad & \dot
q_1=\,\,\,\,c_0q_0-q_0^2+q_1^2+q_2^2+q_3^2,\\
\dot q_2=\,\,\,\,(2q_0-c_0)q_3,& \dot q_3=-(2q_0-c_0)q_2.
\end{array}
\end{equation}
Now
\[ P:=\{q_0=c_0/2\},\quad H=\frac{q_0^2+q_1^2+q_2^2+q_3^2}{
2q_0-c_0},\quad F=q_2^2+q_3^2.
\]
Recall that $P$ is an invariant hyperplane and $H$ and $F$ are two
functionally independent first integrals.

\smallskip

In the invariant space $\mathbb R^4\setminus\{P\}$, we have
$2q_0-c_0\ne 0$. For all regular values $(h,f)$ with $f>0$ and
either $h>c_0$ or $h<0$, we will prove that the invariant torus
$M_{h,f}=\{H=h\}\cap\{F=f\}$ is full of periodic orbits. Taking
the change of coordinates $z=q_0-c_0/2+i\,q_1$, $q_2=r\cos\theta$
and $q_3=r\sin\theta$, equations \eqref{en3.2} become
\begin{equation}\label{3.3}
\dot z=-i\, z^2+i\, \left(r^2+\frac{c_0^2}{4}\right),\quad \dot
r=0,\quad \dot \theta=-2\, \mbox{Re}(z),
\end{equation}
where $\mbox{Re}(z)$ denotes the real part of $z$.
Equations \eqref{3.3} have the solutions
\begin{eqnarray*}
r(t)&=&r,\\
z(t)&=&\frac{\left(z_0+R+(z_0-R) \exp\left(-2Rt\,i\right)\right)
R} {z_0+R-(z_0-R)
\exp\left(-2Rt\,i\right)},\\
\theta(t)&=&\theta_0+2\,\mbox{Re}\left(\int_0^tz(s)ds\right)\\
&=&\theta_0+2Rt+\mbox{Re}\left(\frac 1i\ln\frac{z_0+R-(z_0-R)\exp(-2Rti)}{2R}\right)\\
&=&\theta_0+2Rt+\mbox{Arg}\left(\frac{z_0+R-(z_0-R)\exp(-2Rti)}{2R}\right),
\end{eqnarray*}
with $r\in (0,\infty)$ and $R=\sqrt{r^2+c_0^2/4}$. Clearly, $z(t)$ is a periodic function of period
$\pi/R$ in $t$. Moreover, the third part in the summation of the last equality of $\theta(t)$ is also a periodic
function of period $\pi/R$ in $t$. These show that $q_2$ and $q_3$ are periodic
functions of period $\pi/R$ in $t$, and consequently the orbits on the invariant tori are all periodic.
As a by product of the above proof, we get that with the expansion of the tori their periods become smaller and smaller. This
proves statement $(a)$.

\smallskip

\noindent{\it Statement $(b)$}. Equation \eqref{d4} with $n=3$ is
\begin{equation}\label{d6}
\begin{array}{ll}
\dot q_0= -\left(c_0-3q_0^2+q_1^2+q_2^2+q_3^2\right)q_1,\,\, &
\dot q_1= \left(c_0-q_0^2+3q_1^2+3q_2^2+3q_3^2\right)q_0,\\
\dot q_2=-\left(c_0-3q_0^2+q_1^2+q_2^2+q_3^2\right)q_3, & \dot
q_3= \left(c_0-3q_0^2+q_1^2+q_2^2+q_3^2\right)q_2.
\end{array}
\end{equation}
Now
\[
P:=\{2(q_0^2-q_1^2-q_2^2-q_3^2)-c_0\},
\]
is an invariant generalized hyperboloid, and
\[
H=\frac{(q_0^2+q_1^2+q_2^2+q_3^2)^2}{2(q_0^2-q_1^2-q_2^2-q_3^2)-c_0},\quad
F=q_2^2+q_3^2,
\]
are two functionally independent first integrals.

\smallskip

For each regular values $(h,f)$ of $(H,F)$, we study the dynamics
on the invariant tori $M_{h,f}$. Since $f>0$, the generalized
cylinder $F=f$ is parametrized by
\[
q_2=\sqrt{f}\cos\theta,\quad q_3=\sqrt{f}\sin\theta.
\]
Restricted to the $F=f$, the hypersurface $H=h$ with
$f^2+(2f+c_0)h<0$ can be parametrized by
\[
q_0=\sqrt{G(\v,h)}\cos \v,\quad q_1=\sqrt{G(\v,h)}\sin\v,
\]
with
\[
G(\v,h)=h\cos 2\v-f+\sqrt{h^2\cos^22\v-2fh\cos 2\v-2fh-c_0h}.
\]
Note that here we study only those tori $M_{h,f}$ with $f>0$ and
$f^2+(2f+c_0)h<0$. They are probably the most simple ones which
can be parametrized.

\smallskip

On the above mentioned invariant torus $M_{h,f}$, system
\eqref{d6} writes in
\begin{equation}\label{d7}
\begin{array}{l}
\dot h=0,\quad  \dot\theta=c_0+f-G(\v,h)(2\cos 2\v+1),\\
\dot r=0,\quad  \dot\v=c_0-G(\v,h)\cos2\v+f\cos2\v+2f.
\end{array}
\end{equation}
Set
\begin{eqnarray*}
A(\v,h)&=&\sqrt{h^2\cos^22\v-2fh\cos 2\v-2fh-c_0h},\\
B(\v,h)&=&c_0+f-G(\v,h)(2\cos 2\v+1).
\end{eqnarray*}
Then
\[
c_0-G(\v,h)\cos2\v+f\cos2\v+2f=\frac{AB}{A+2h\cos^2\v}.
\]
If there is a periodic orbit on $M_{h,f}$, we assume that its
smallest positive period is $2m\pi$ in $\v$ and $2n\pi$ in
$\theta$ for $m, n\in\mathbb N$. We get from \eqref{d7} that
\[
\int_0^{2m\pi}\left(1+\frac{2h\cos^2\v}{A}\right)d\v=\int_0^{2n\pi}d\theta.
\]
The last equality can be written in
\[
2(n-m)\pi=m\int_0^{2\pi}\frac{h(1+\cos^2\psi)}{\sqrt{h^2\cos^2\psi-2fh\cos\psi-(2f+c_0)h}}d\psi.
\]
We can check easily that
\[
I(h):=\frac
1{2\pi}\int_0^{2\pi}\frac{h(1+\cos^2\psi)}{\sqrt{h^2\cos^2\psi-2fh\cos\psi-(2f+c_0)h}}d\psi,
\]
is analytic in $h$ with $f^2+(2f+c_0)h<0$, and $I'(h)\not\equiv 0$ in any open subset of $\mathbb R$. The last claim
implies that $I(h)$ is a locally open mapping. So, for any given $f>0$ there exist infinitely many
$h$ such that $M_{h,f}$ is full of periodic orbits, and also
infinitely many  $h$ for which $M_{h,f}$ has dense orbits. This
proves statement $(b)$.

\smallskip

We complete the proof of the theorem. $\qquad \qquad \Box$

\bigskip

\subsection{Proof of Theorem \ref{th2.5}}

Working in a similar way to the proof of Theorem \ref{t1}, we only
need to study equation \eqref{e5} with $c=c_0+ic_1$. Writing
equation \eqref{e5} in a system gives
\begin{equation}\label{e50}
\begin{array}{ll}
\dot q_0=c_0q_0-c_1q_1-q_0^2+q_1^2+q_2^2+q_3^2,\,\, & \dot
q_1=c_1q_0+(c_0-2q_0)q_1,\\
\dot q_2=(c_0-2q_0)q_2-c_1q_3, & \dot q_3=c_1q_2+(c_0-2q_0)q_3.
\end{array}
\end{equation}
Recall that $L=c_0 q_0+c_1q_1-(c_0^2+c_1^2)/2$.

\smallskip

\noindent $(a)$ Restricted to the hyperplane $L=0$ we have
\[
\left.\frac{dL}{dt}\right|_{\eqref{e50}}=4c_0(q_0^2+q_1^2+q_2^2+q_3^2).
\]
So every orbit intersects the hyperplane $L=0$ transversally.
Obviously $L=0$ is orthogonal to the line connecting the
singularities $O$ and $S$, and have the same distance to $O$ and
$S$.

\smallskip

Set
\begin{equation}\label{e51}
H=\frac{q_0^2+q_1^2+q_2^2+q_3^2}{2c_0q_0+2c_1q_1-K_0},
\end{equation}
where $K_0=c_0^2+c_1^2$. We have
\begin{equation}\label{e52}
\left.\frac{dH}{dt}\right|_{\eqref{e50}}=\frac{-2c_0(q_0^2+q_1^2+q_2^2+q_3^2)
B} {\left(2c_0q_0+2c_1q_1-K_0\right)^2},
\end{equation}
where $B=(q_0-c_0)^2+(q_1-c_1)^2+q_2^2+q_3^2$. The level set $H=h$
is empty if $0<h<1$, and is a ball, denoted by $B_h$, centered at
$(hc_0,hc_1,0,0)$ with the radius $\sqrt{(c_0^2+c_1^2)(h^2-h)}$ if
$h>1$ or $h<0$. We can check easily that the balls $B_h$ with
$h>1$ (resp. $h<0$) contain the singularity $S$ (resp. $O$) in
their interiors and are located in $L>0$ (resp. $L<0$).
Furthermore, it is easy to prove that when $h\searrow 1$ (resp.
$h\nearrow 0$) the ball $B_h$ shrinks to the singularity $S$
(resp. $O$), and that when $h\nearrow \infty$ (resp. $h\searrow
-\infty$) the ball $B_h$ expands and approaches the hyperplane
$L=0$.

\smallskip

From the derivative of $H$ and the property of the ball $B_h$, it
follows that each orbit starting on $L=0$ will be heteroclinic
connecting the two singularities $S$ and $O$. Moreover we get from
the last two equations of \eqref{e50} that these orbits spirally approach
$S$ and $O$. These last proofs imply that except those
orbits being heteroclinic to $S$ and $O$, there are two other
ones: one is heteroclinic to $S$ and infinity, and another is
heteroclinic to $O$ and infinity. This proves statement $(a)$.

\smallskip

\noindent {\it Statement $(b)$}. Since $c_0=0$, we get from
\eqref{e52} that the function $H$ defined in \eqref{e51} is a
first integral of system \eqref{e50}. Furthermore we can prove
that
\[
F=\frac{(q_0^2+q_1^2+q_2^2+q_3^2)^2}{(2q_1-c_1)^2+4q_2^2+4q_3^2},
\]
is also a first integral of system \eqref{e50}. In addition,
$L=0$, i.e. $2q_1=c_1$ is invariant.

\smallskip

Some calculations show that $H$ and $F$ are functionally
independent, and that the critical points are
$\{q_2=0\}\cap\{q_3=0\}$ and
$\{q_0=0\}\cap\{-a_1q_1+q_1^2+q_2^2+q_3^2=0\}$. The corresponding
critical values of $(H,F)$ are $(h,0)$ and $(h,f)$ with $f>0$ and
$h=c_1\sqrt{f}/(2\sqrt f-c_1)$. On the invariant plane
$\{q_2=0\}\cap\{q_3=0\}$, there are two period annuli separated by
the invariant line $L=0$. On the invariant sphere
$\{q_0=0\}\cap\{-a_1q_1+q_1^2+q_2^2+q_3^2=0\}$ all orbits are
periodic. This last claim follows from the fact that system
\eqref{e50} restricted to the sphere has the first integral
$q_2^2+q_3^2$.

\smallskip

For any $h>1$ or $h<0$ and $f>0$ with $h\ne c_1\sqrt{f}/(2\sqrt
f-c_1)$, the values $(h,f)$ are regular for $(H,F)$. Since the
hypersurfaces $H=h$ and $F=f$ are compact and intersect
transversally, their intersections denoted by $M_{h,f}$ should be
two dimensional compact invariant manifolds (if exists). We claim
that the connected parts of $M_{h,f}$ are $2$--dimensional
invariant tori.

\smallskip

We now prove the claim. For doing so, it suffices to show that
$M_{h,f}$ is orientable and has genus $1$. Associated to the
$2$--field
\[
\begin{array}{l}
\nabla H\wedge\nabla F=\left(\frac{\partial H}{\partial
q_0}\frac{\partial F}{\partial q_1}-\frac{\partial H}{\partial
q_1}\frac{\partial F}{\partial q_0}\right)\frac{\partial
}{\partial q_0}\wedge\frac{\partial }{\partial q_1}+
\left(\frac{\partial H}{\partial q_0}\frac{\partial F}{\partial
q_2}-\frac{\partial H}{\partial q_2}\frac{\partial F}{\partial
q_0}\right)\frac{\partial }{\partial q_0}\wedge\frac{\partial
}{\partial q_2}\\
\qquad \qquad\quad+\left(\frac{\partial H}{\partial
q_0}\frac{\partial F}{\partial q_3}-\frac{\partial H}{\partial
q_3}\frac{\partial F}{\partial q_0}\right)\frac{\partial
}{\partial q_0}\wedge\frac{\partial }{\partial
q_3}+\left(\frac{\partial H}{\partial q_1}\frac{\partial
F}{\partial q_2}-\frac{\partial H}{\partial q_2}\frac{\partial
F}{\partial q_1}\right)\frac{\partial
}{\partial q_1}\wedge\frac{\partial }{\partial q_2}\\
\qquad \qquad\quad+\left(\frac{\partial H}{\partial
q_1}\frac{\partial F}{\partial q_3}-\frac{\partial H}{\partial
q_3}\frac{\partial F}{\partial q_1}\right)\frac{\partial
}{\partial q_1}\wedge\frac{\partial }{\partial
q_3}+\left(\frac{\partial H}{\partial q_2}\frac{\partial
F}{\partial q_3}-\frac{\partial H}{\partial q_3}\frac{\partial
F}{\partial q_2}\right)\frac{\partial }{\partial
q_2}\wedge\frac{\partial }{\partial q_3},
\end{array}
\]
the dual $2$--form is
\[
\begin{array}{l}
\omega=\quad\left(\frac{\partial H}{\partial q_2}\frac{\partial
F}{\partial q_3}-\frac{\partial H}{\partial q_3}\frac{\partial
F}{\partial q_2}\right)dq_0dq_1- \left(\frac{\partial H}{\partial
q_1}\frac{\partial F}{\partial q_3}-\frac{\partial H}{\partial
q_3}\frac{\partial F}{\partial
q_1}\right)dq_0dq_2\\
\quad \quad+\left(\frac{\partial H}{\partial q_1}\frac{\partial
F}{\partial q_2}-\frac{\partial H}{\partial q_2}\frac{\partial
F}{\partial q_1}\right)dq_0dq_3 +\left(\frac{\partial H}{\partial
q_0}\frac{\partial F}{\partial q_3}-\frac{\partial H}{\partial
q_3}\frac{\partial F}{\partial
q_0}\right)dq_1dq_2\\
\quad\quad- \left(\frac{\partial H}{\partial q_0}\frac{\partial
F}{\partial q_2}-\frac{\partial H}{\partial q_2}\frac{\partial
F}{\partial q_0}\right)dq_1dq_3+
 \left(\frac{\partial H}{\partial
q_0}\frac{\partial F}{\partial q_1}-\frac{\partial H}{\partial
q_1}\frac{\partial F}{\partial q_0}\right)dq_2dq_3.
\end{array}
\]
Recall that $\nabla$ denotes the gradient of a smooth function.
Since the fields $\nabla H$ and $\nabla F$ are linearly
independent on $M_{h,f}$, the two form $\omega$ is non--zero on
$M_{h,f}$. Hence $M_{h,f}$ is orientable, see e.g. \cite[Sec.
2.5]{AM}  and also \cite{CGM}.

\smallskip

Denote by $\mathcal X_{h,f}$ the restriction of the vector field
defined by \eqref{d4} to $M_{h,f}$. Since the vector field
$\X_{h,f}$ has no singularities, applying the Poincar\'e--Hopf
formula to the manifold $M_{h,f}$ we get
\[
0=\mbox{ind}(\mathcal X_{h,h_1})=\chi(M_{h,h_1})=2-2g,
\]
where  ind$(\mathcal X_{h,f})$ denotes the sum of the indices of
the singularities of $\mathcal X_{h,f}$ on $M_{h,f}$, and
$\chi(M_{h,f})$ and $g$ are the Euler characteristic and the genus
of the surface $M_{h,f}$, respectively. This shows that the genus
of $M_{h,f}$ is one. It is well--known that an orientable compact
connected surface of genus one is a torus, see e.g., \cite[Sec.
X]{La} for more details.

\smallskip

We complete the proof of statement $(b)$ and consequently the
proof of the theorem.$\qquad \qquad \Box$

\bigskip

\subsection{Proof of Theorem \ref{th2.6}}

\smallskip
Recall that $L=q_0^2-q_1^2-q_2^2-q_3^2-c_0^2/2$. For simplifying
the notations we denote by $H_+$ and $H_-$ the subset of $\R^4$
with $L>0$ and $L<0$ respectively, and by $H_+^+$ and $H_+^-$ the
two parts of $H_+$ with $q_0>c_0/\sqrt 2$ and $q_0<-c_0/\sqrt 2$,
respectively.

\smallskip

Working in a similar way to the proof of Theorem \ref{t1}, we
assume without loss of generality that $a=a_0+a_1i$. Equation
\eqref{e3} is equivalent to the system
\begin{equation}\label{e4.15}
\begin{array}{ll}
\dot q_0=\,\,\,\,\,a_1q_1A -a_0q_0B,\qquad & \dot q_1=-a_0q_1A
-a_1q_0B,\\
\dot q_2=-(a_0q_2-a_1q_3)A, & \dot q_3=-(a_1q_2+a_0q_3)A ,
\end{array}
\end{equation}
where $A=c_0^2-3q_0^2+q_1^2+q_2^2+q_3^2$ and
$B=c_0^2-q_0^2+3(q_1^2+q_2^2+q_3^2)$. It is easy to check that
system \eqref{e4.15} has the three finite singularities $O$, $S_+$
and $S_-$. Furthermore restricted to $L=0$ the derivative of $L$
along the solutions of system \eqref{e4.15} with respect to the
time $t$ is
\begin{equation}\label{e4.16}
\left.\frac{dL}{dt}\right|_{\eqref{e4.15}}=-2a_0(2q_0^2-c_0^2/2)^2.
\end{equation}

\smallskip

\noindent {\it Statement $(a)$}. We consider the case $a_0<0$. The
proof of the case $a_0>0$ follows from the same arguments than that of $a_0<0$. By
\eqref{e4.16} we get that if an orbit intersects $L=0$, it should
transversally pass through it. Moreover the orbits meeting $L=0$ will go
from $H_-$ to $H_+$ as the time increases.

\smallskip

Set
\begin{equation}\label{e4.17}
\begin{array}{l}
H=(q_0^2+q_1^2+q_2^2+q_3^2)^2/(q_0^2-q_1^2-q_2^2-q_3^2-c_0^2/2).
\end{array}
\end{equation}
We can check that for $h\in(2c_0^2,\infty)$ the hypersurface
$E_h:=\{H=h\}$ has two branches which are located in $H_+^+$ and
$H_+^-$ respectively, and that for $h\in(-\infty,0)$ the
hypersurface $E_h$ has a unique branch which is located in $H_-$.
Moreover we can check that for
$h\in(2c_0^2,\infty)\cup(-\infty,0)$ the hypersurface $E_h$ is
compact and contains one of the three singularities in its
interior. When $h\rightarrow \pm\infty$ the hypersurface $E_h$
approaches the hyperboloid $L=0$.

\smallskip

Now we can verify that
\begin{equation}\label{e4.18}
\left.\frac{dH}{dt}\right|_{\eqref{e4.15}}=8a_0N\frac{(q_0^2+q_1^2+q_2^2+q_3^2)^2}
{\left(c_0^2-2(q_0^2-q_1^2-q_2^2-q_3^2)\right)^2},
\end{equation}
where $
N=\left(c_0^4-2c_0^2(q_0^2-q_1^2-q_2^2-q_3^2)+(q_0^2+q_1^2+q_2^2+q_3^2)^2\right)
$. Since outside $S_1$ and $S_2$ we have $N>0$,  it follows that
the subsets $H_+$ (resp. $H_-$) are positively (resp. negatively)
invariant by the flow of the system. Furthermore, all orbits
starting in $H_+^+$ (resp. $H_+^-$) will approach $S_+$ (resp.
$S_-$) when $t\rightarrow \infty$. All orbits starting in $H_-$
will go to $O$ when $t\rightarrow -\infty$. So all orbits starting
on $L=0$ with $q_0>c_0/\sqrt 2$ (resp. $q_0<-c_0/\sqrt 2$) will be
heteroclinic to $O$ and $S_+$ (resp. to $O$ and $S_-$). This
proves statement $(a)$.

\smallskip

\noindent {\it Statement $(b)$}. Equations \eqref{e4.16} and
\eqref{e4.18} show that the hyperboloid $L=0$ is invariant and
that $H$ is a first integral of system \eqref{e4.15}. Moreover we
can prove that
\[
F=q_2^2+q_3^2,
\]
is also a first integral of \eqref{e4.15}, and that $H$ and $F$
are functionally independent. Recall that $F=f$ is a
$3$--dimensional cylinder when $f>0$ and is a plane when $f=0$.

\smallskip

Working in a similar way to the proof of statement $(b)$ of
Theorem \ref{th2.5} we can prove that for
$h\in(-\infty,0)\cup(2c_0^2,\infty)$ the intersections
$E_h\cap\{F=f\}$ are either formed by periodic orbits for $(h,f)$
being critical values or $2$--dimensional invariant tori for $(h,
f)$ being regular values. Using the same methods as those given in
the proof of statement $(b)$ of Theorem \ref{t1.1} we can prove
that of the invariant tori there are infinitely many ones
fulfilling periodic orbits and also infinitely many ones
fulfilling dense orbits. This proves statement $(b)$ and
consequently the theorem. \qquad $\Box$

\bigskip

\section{Appendix: the linear case}

\smallskip

For the homogeneous linear differential equations
\begin{equation}\label{esu1}
\dot q=a q+q b,
\end{equation}
with $a,b\in\mathbb H$ nonzero, taking $H=q\overline q$ we have
\[
\left.\frac{dH}{dt}\right|_{\eqref{esu1}}=(a+\overline
{a}+b+\overline b)(q\overline q).
\]
Moreover its equivalent $4$--dimensional linear differential system
has at the origin the four eigenvalues
\[
(a+\overline a+b+\overline b)/{2}\pm\left(\sqrt{(a-\overline
a)^2}\pm\sqrt{(b-\overline b)^2}\right)/{2}.\] So the dynamics of
\eqref{esu1} follows easily from these eigenvalues.

\smallskip

For the homogeneous linear differential equations
\begin{equation}\label{esu2}
\dot q=a q+\overline q b,
\end{equation}
with $a,b\in\mathbb H$ nonzero, its equivalent $4$--dimensional
linear differential system has the four eigenvalues
\[
(a-b+\overline{ a-b})/{2}\pm\sqrt{(a+b-\overline {a+b})^2}\,/2,\quad
(a+\overline a)/{2}\pm\sqrt{(a-\overline a)^2-b\overline b}\,/{2}.
\]
Then the dynamics of \eqref{esu2} follows easily from these
eigenvalues.

\smallskip

For the homogeneous linear equations
\begin{equation}\label{esu3}
\dot q=a q+ b \overline q,
\end{equation}
with $a,b\in\mathbb H$ nonzero, its equivalent $4$--dimensional
linear differential system has the four eigenvalues
\[
(a-b+\overline{ a-b})/{2}\pm\sqrt{(a-b-\overline {a-b})^2}\,/2,
\quad (a+\overline a)/{2}\pm\sqrt{(a-\overline a)^2-b\overline
b}\,/{2}.
\]
Then its dynamics follows also from these eigenvalues.

\smallskip

For the non--homogeneous linear quaternion differential equations
\begin{equation}\label{e2.1}
\dot q=b+aq, \qquad \dot q=b+qa,
\end{equation}
with $a,b\in \mathbb H$ nonzero, they can be transformed to
homogeneous ones via the change of variables $p=q+a^{-1}b$ or
$p=q+ba^{-1}$. So their dynamics can be obtained from Theorem 2 of
\cite{GLZ}.

\smallskip

\smallskip

\subsection*{Acknowledgements}

The author thanks Professors Armengol Gasull and Jaume Llibre for
their discussion and comments to part of results given in the
first version of this paper. I should appreciate the referees for their
excellent comments and suggestions, which can improve our paper both in mathematics and in the expressions.


\begin{thebibliography}{99}


\bibitem{AM} {\sc R. Abraham and J.E. Marsden},  {\it Foundations of Mechanics} 2nd Ed.,
Addison--Wesley, Redwood City, California, 1987.

\bibitem{Ad1} {\sc S.L. Adler}, {\it Quaternionic quantum field theory},
Commun. Math. Phys. {\bf 104} (1986), 611--656.

\bibitem{Ad2} {\sc S.L. Adler}, {\it Quaternionic Quantum Mechanics and Quantum
Fields},  Oxford University Press, New York, 1995.

\bibitem{AGP} {\sc M.J.\'{A}lvarez, A. Gasull and R. Prohens}, {\it Configurations of critical
points in complex polynomial differential equations}, Nonlinear
Analysis {\bf 71} (2009), 923--934.


\bibitem{Ar} {\sc V.I. Arnold},  {\it Mathematial Methods of Classical
Mechanics}, Springer-Verlag, New York, 1978.

\bibitem{BC04} {\sc M. Bruschi and F. Calogero}, {\it Integrable systems of quartic oscillators. II},
Phys. Lett. A {\bf 327} (2004), 320--326.

\bibitem{CD06} {\sc F. Calogero and A. Degasperis}, {\it New integrable PDEs of boomeronic type},
J. Phys. A {\bf 39} (2006),8349--8376.

\bibitem{CD04} {\sc F. Calogero and A. Degasperis}, {\it New integrable equations of nonlinear Schr\"{o}dinger type},
Stud. Appl. Math. {\bf 113} (2004), 91--137.

\bibitem{CCZ} {\sc C. Chen, J. Cao and X. Zhang}, {\it The topological
structure of the Rabinovich system having an invariant algebraic
surface}, Nonlinearity {\bf 21} (2008), 211--220.


\bibitem {CM}
{\sc J. Campos and J. Mawhin}, {\it Periodic solutions of
quaternionic-values ordinary differential equations}, Annali di
Matematica {\bf 185} (2006), S109--S127.

\bibitem {CGM}
{\sc A. Cima, A. Gasull and V. Ma\~{n}osa}, {\it Some properties
of the $k$--dimensional Lyness's map}, J. Phys. A: Math. Theor.
{\bf 41} (2008), 285205.

\bibitem{FJSS} {\sc D. Finkelstein, J.M. Jauch, S. Schiminovich and D. Speiser},
{\it Foundations of quaternion quantum mechanics}, J. Math. Phys.
{\bf 3} (1962), 207--220.

\bibitem {Fr} {\sc F.G. Frobenius},  {\it Ueber lineare
Substitutionen und bilineare Formen}, J. Reine Angew. Math. {\bf
84} (1878), 1-63.

\bibitem{GLMM} {\sc A. Gasull, J. Llibre, V. M\~{a}osa and F. M\~{a}osas},
{\it The focus-centre problem for a type of degenerate system},
Nonlinearity {\bf 13} (2000),  699--729.


\bibitem{GLZ} {\sc A. Gasull. J. Llibre and Xiang Zhang}, {\it One--dimensional quaternion homogeneous polynomial
differential equations}, J. Math. Phys. {\bf 50} (2009), 082705.


\bibitem {Gi}
{\sc J.D. Gibbon}, {\it A quaternionic structure in the
three--dimensional Euler and ideal magneto--hydrodynamics
equation}, Physica D {\bf 166} (2002), 17--28.

\bibitem {GHKR}
{\sc J.D. Gibbon, D.D. Holm, R.M. Kerr and I. Roulstone}, {\it
Quaternions and particle dynamics in the Euler fluid equations},
Nonlinearity {\bf 19} (2006), 1969--1983.


\bibitem{GL04} {\sc J. Gin\'e and J. Llibre}, {\it
Integrability and algebraic limit cycles for polynomial
differential systems with homogeneous nonlinearities}, J.
Differential Equations {\bf 197} (2004), 147--161.

\bibitem {Ha} {\sc S.W.R. Hamilton}, {\it Lectures on Quaternions}, Royal
Irish Academy, Hodges and Smith, Dublin, 1853.

\bibitem{Han} {\sc A.J. Hanson}, {\it Visualizing Quaternions}, Elsevier,
San Francisco, 2006.

\bibitem{IC02} {\sc S. Iona and F. Calogero}, {\it Integrable systems of quartic oscillators in ordinary (three-dimensional) space},
J. Phys. A {\bf 35} (2002), 3091--3098.

\bibitem{La} {\sc S. Lang},  {\it Differential and Riemannian Manifolds},
Springer-Verlag, New York, 1995.

\bibitem{LLLZ} {\sc C. Li, W. Li, J. Llibre and Z. Zhang}, {\it On the
limit cycles of polynomial differential systems with homogeneous
nonlinearities},  Proc. Edinburgh Math. Soc. (2) {\bf 43} (2000),
529--543.

\bibitem{Ll} {\sc J. Llibre}, {\it Handbook of Differential Equations}, Elsevier/North--Holland,
Amsterdam, 2004, pp. 437--532.

\bibitem{LV09} {\sc J. Llibre and C. Valls}, {\it Classification of the centers, their
cyclicity and isochronicity for a class of polynomial differential
systems generalizing the linear systems with cubic homogeneous
nonlinearities}, J. Differential Equations {\bf 246} (2009),
2192--2204.

\bibitem{LY} {\sc J. Llibre and J. Yu}, {\it On the periodic orbits of the static,
spherically symmetric Einstein-Yang-Mills equations}, Commun.
Math. Phys. {\bf 286} (2009), 277--281.

\bibitem{LZ1} {\sc J. Llibre and Xiang Zhang}, {\it Invariant algebraic surfaces of the Lorenz systems},
J. Math. Phys.  {\bf 43} (2002), 1622-1645..


\bibitem {LZ}
{\sc J. Llibre and Xiang Zhang}, {\it Darboux theory of
integrability in $C^n$ taking into account the multiplicity}, J.
Differential Equations {\bf 246} (2009), 541--551.

\bibitem {RR97}
{\sc V.N. Roubtsov and I. Roulstone}, {\it Examples of
quaternionic and K$\ddot{\mbox{a}}$hler structures in Hamiltonian
models of nearly geostrophic flow}, J. Phys. A: Math. Gen. {\bf
30} (1997), L63--L68.

\bibitem {RR01}
{\sc V.N. Roubtsov and I. Roulstone}, {\it Holomorphic structures in
hydrodynamical models of nearly geostrophic flow}, Proc. R. Soc.
London A {\bf 457} (2001), 1519--1531.

\bibitem {Wi09}
{\sc P. Wilczynski}, {\it Quaternionic--valued ordinary
differential equations. The Riccati equation}, J. Differential
Equations {\bf 247} (2009), 2163--2187.

\end{thebibliography}
\end{document}